# **Demystification of Taylor , Laurent Coefficients of Lerch , Hurwitz Zeta Functions and Dirichlet L-Function at Unity and Zero and their Bounds**

**VIVEK VISHWANATH RANE**

Department of Mathematics ,
The Institute of Science,
15 , Madam Cama Road,
Mumbai-400 032.
INDIA
v_v_rane@yahoo.co.in

**Abstract :** Using elementary methods , we obtain simple explicit expressions , bounds and the interrelations of Laurent coefficients of Hurwitz Zeta function , Taylor coefficients of Dirichlet L-series and the Taylor coefficients of Lerch's Zeta function at $s=1$ (and also at $s=0$) . In particular , we obtain simple , explicit expression and very good upper bound for $L^{(r)}(1,\chi)$ as a function of both $q$ and $r$ , where $L(s,\chi)$ is a Dirichlet L-function corresponding to a Dirichlet character $\chi$ modulo $q$ and the superscript $(r)$ denotes $r-th$ order derivative with respect to s . We also give elementary complete forms of approximate functional equations of $\zeta(s,\alpha)$ and $L(s,\chi)$ and their derivatives ,in the critical strip .

**Keywords** : Hurwitz/Lerch Zeta function , Dirichlet character , Dirichlet L-series .

# **Demystification of Taylor , Laurent Coefficients of Lerch , Hurwitz Zeta Functions and Dirichlet L-Function at Unity and Zero and their Bounds**

**VIVEK VISHWANATH RANE**
**v_v_rane@yahoo.co.in**

Introduction **:** Let $s=\sigma+it$ be a complex variable , where $\sigma$ and $t$ are real . For $\alpha$ with $0<\alpha\le 1$ , let $\zeta(s,\alpha)$ be Hurwitz zeta function defined by $\zeta(s,\alpha)=\sum_{n\ge 0}(n+\alpha)^{-s}$ for $\sigma>1$ ; and its analytic continuation . For any integer $q\ge 1$ , let $\chi(\bmod q)$ be a Dirichlet character and let $L(s,\chi)$ be the corresponding Dirichlet L-series . Thus $L(s,\chi)=\sum_{n\ge 1}\chi(n)n^{-s}$ for $\sigma>1$ ; and its analytic continuation . Note that $L(s,\chi)=\sum_{a=1}^{q}\chi(a)Z(s,a,q)$ , where $Z(s,a,q)=q^{-s}\cdot\zeta(s,\frac{a}{q})$ . In what follows , we shall be considering $L(s,\chi)$ for a non-principal Dirichlet character only . For $0\le\lambda\le 1$ and for $0<\alpha\le 1$ , let $\phi(\lambda,\alpha,s)$ be Lerch's zeta function defined by $\phi(\lambda,\alpha,s)=\sum_{n\ge 0}e^{2\pi i\lambda n}(n+\alpha)^{-s}$ for $\sigma>1$ ; and its analytic continuation . Note that if $\lambda$ is an integer , $\phi(\lambda,\alpha,s)=\zeta(s,\alpha)$ , the Hurwitz's zeta function . Also , if $\lambda$ is not an integer , then $\phi(\lambda,\alpha,s)$ is an entire function of $s$ . For any integer $r\ge 0$ , let the superscript $(r)$ denote differentiation with respect to $s$ so that $\zeta^{(r)}(s,\alpha), Z^{(r)}(s,a,q), \phi^{(r)}(\lambda,\alpha,s)$ and $L^{(r)}(s,\chi)$ stand for $r-th$ order derivatives with respect to s of the respective functions . We shall write $[u]$ for the integral part of a real variable $u$ . We write $\psi(u) = u - [u] - \frac{1}{2}$ . Note that $\psi(u) = -\sum_{|n|\geq 1} \frac{e^{2\pi i n u}}{2\pi i n}$ for non-integral $u$ . We write $\psi_2(u) = -\sum_{|n|\geq 1} \frac{e^{2\pi i n u}}{(2\pi i n)^2}$ so that $\frac{d}{du}\psi_2(u) = \psi(u)$ . Incidentally note that $(-1)^r \phi^{(r)}(\lambda,\alpha,s)$ $= \sum_{n\geq 0} e^{2\pi i \lambda n}(n+\alpha)^{-s} \log^r(n+\alpha)$ for $\sigma > 1$ ; and its analytic continuation . In what follows , we shall use the fact that $\sum_{a\leq x} \chi(a) << q^{1/2} \log q$ for any non-principal Dirichlet character $\chi(\bmod q)$ and for any positive x . In what follows , $\Gamma(s)$ stands for gamma function and the logarithm stands for the principal value of logarithm . Also , O- and << - constants will be absolute , unless stated otherwise .

It is known that $\zeta(s,\alpha) - \frac{1}{s-1} = \sum_{r\geq 0} \gamma_r(\alpha)(s-1)^r$ , where the coefficients $\gamma_r(\alpha)$ are given by $\gamma_r(\alpha) = \frac{(-1)^r}{r!} \lim_{N\to\infty} \left\{ \sum_{n=0}^{N} \frac{\log^r(m+\alpha)}{n+\alpha} - \frac{\log^{r+1}(N+\alpha)}{r+1} \right\}$ .

Evidently , $\gamma_r(\alpha) = \frac{1}{r!} \lim_{s\to 1} \left( \zeta(s,\alpha) - \frac{1}{s-1} \right)^{(r)}$ .

$$= \frac{1}{r!}\lim_{s\to 1}\left(\zeta^{(r)}(s,\alpha) - \left(\frac{1}{s-1}\right)^{(r)}\right) = \frac{1}{r!}\lim_{s\to 1}\left(\zeta^{(r)}(s,\alpha) - \frac{(-1)^r r!}{(s-1)^{r+1}}\right).$$

On the other hand , it is known that $L^{(r)}(1,\chi) = \sum_{a=1}^{q}\chi(a)\cdot(-1)^r\gamma_r(a,q)$, where

$$\gamma_r(a,q) = \lim_{N\to\infty}\left\{\sum_{\substack{n=a(\mathrm{mod}\,q)\\ n\le N}}\frac{\log^r n}{n} - \frac{\log^{r+1} N}{q(r+1)}\right\}$$ for any non-principal character $\chi(\mathrm{mod}\,q)$ .

This means , if $\chi(\mathrm{mod}\,q)$ is a non-principal Dirichlet character , then

$L(s,\chi) = \sum_{r\ge 0}\gamma_r(\chi)(s-1)^r$ , where $\gamma_r(\chi) = \frac{L^{(r)}(1,\chi)}{r!} = \sum_{a=1}^{q}\chi(a)\left(\frac{(-1)^r}{r!}\gamma_r(a,q)\right)$ .

Actually , $L(s,\chi) = \sum_{a=1}^{q}\chi(a)Z(s,a,q)$ , where $Z(s,a,q) = q^{-s}\cdot\zeta(s,\frac{a}{q})$ so that $L^{(r)}(s,\chi) = \sum_{a=1}^{q}\chi(a)\cdot Z^{(r)}(s,a,q)$ .

In particular , $L^{(r)}(1,\chi) = \sum_{a=1}^{q}\chi(a)\cdot\lim_{s\to 1} Z^{(r)}(s,a,q) = \sum_{a=1}^{q}\chi(a)\cdot\lim_{s\to 1}\left(q^{-s}\zeta(s,\frac{a}{q})\right)^{(r)}$

$$= \sum_{a=1}^{q}\chi(a)\lim_{s\to 1}\left\{(q^{-s}(\zeta(s,\tfrac{a}{q}) - \tfrac{1}{s-1})^{(r)} + (\tfrac{q^{-s}}{s-1})^{(r)}\right\}$$

$$= \sum_{a=1}^{q}\chi(a)\lim_{s\to 1}\left\{q^{-s}(\zeta(s,\tfrac{a}{q}) - \tfrac{1}{s-1})\right\}^{(r)} + \lim_{s\to 1}\sum_{a=1}^{q}\chi(a)\left(\tfrac{q^{-s}}{s-1}\right)^{(r)}$$

$$= \sum_{a=1}^{q}\chi(a)\lim_{s\to 1}\left\{q^{-s}(\zeta(s,\tfrac{a}{q}) - \tfrac{1}{s-1})\right\}^{(r)} = \sum_{a=1}^{q}\chi(a)\cdot(-1)^r\gamma_r(a,q) .$$

In fact , it can be seen that $(-1)^r\gamma_r(a,q) = \lim_{s\to 1}\left\{q^{-s}\left(\zeta(s,\frac{a}{q}) - \frac{1}{s-1}\right)\right\}^{(r)}$

The exact connection between $\gamma_r(a,q)$ and $\gamma_r(\frac{a}{q})$ can be seen from the following

Proposition , which will be proved by two methods .

**Proposition** : Let $\zeta(s,\alpha) = \frac{1}{s-1} + \sum_{r\geq 0} \gamma_r(\alpha)(s-1)^r$ and let $L(1,\chi) = \sum_{a=1}^{q} \chi(a) \cdot (-1)^r \gamma_r(a,q)$

with $\chi(\bmod q)$, a non-principal character modulo an integer $q \geq 1$, where

$(-1)^r \gamma_r(a,q) = \lim_{s\to 1} \left\{ q^{-s}\left(\zeta(s,\frac{a}{q}) - \frac{1}{s-1}\right) \right\}^{(r)}$. Then

$$(-1)^r \gamma_r(a,q) = \frac{r!}{q} \sum_{\ell=0}^{r} \frac{(-1)^\ell}{\ell!} \log^\ell q \cdot \gamma_{r-\ell}(\tfrac{a}{q}) .$$

Thus the problem boils down to finding simple explicit expressions for the coefficients of Laurent series of $\zeta(s,\alpha)$ at $s=1$; and the coefficients of Taylor expansion of $L(s,\chi)$ at $s=1$ or equivalently finding simple , explicit expressions for $L^{(r)}(1,\chi)$ for $r \geq 0$; and estimating them . There is nothing sacrosanct of Taylor/Laurent coefficients at $s=1$ . We may , as well consider the Taylor coefficients of $\zeta(s,\alpha)$ in the form $\zeta(s,\alpha) = \sum_{r\geq 0} \beta_r(\alpha) s^r$ or $L(s,\chi) = \sum_{r\geq 0} \beta_r(\chi) s^r$ . More generally for $0 < \lambda < 1$ and for $0 < \alpha \leq 1$ , we may consider Taylor series expressions for $\phi(\lambda,\alpha,s)$ in the form $\phi(\lambda,\alpha,s) = \sum_{r\geq 0} \gamma_r(\lambda,\alpha)(s-1)^r$ or $\phi(\lambda,\alpha,s) = \sum_{r\geq 0} \beta_r(\lambda,\alpha) s^r$ . This amounts to dealing with $\phi^{(r)}(\lambda,\alpha,s)$ for $s=0$ or $1$ . Our Theorem 1 below shows that in the critical strip $0<\sigma\leq 1$, expressions for $\zeta^{(r)}(s,\alpha)$, $(q^{-s}\zeta(s,\alpha))^{(r)}$ , $L^{(r)}(s,\chi)$ or $\phi^{(r)}(\lambda,\alpha,s)$ can be obtained from the expressions of $\zeta(s,\alpha)$ , $q^{-s}\zeta(s,\frac{a}{q})$ , $L(s,\chi)$ and $\phi(\lambda,\alpha,s)$ respectively , by differentiating $r$ times with respect to $s$ , term-by-term and by differentiating under integral sign , if need arises .

There are some results giving upper bounds for $L^{(r)}(1,\chi)$ ( with non-principal Dirichlet character $\chi(\mathrm{mod}\,q)$ ) , when q and $\chi(\mathrm{mod}\,q)$ are of certain types only . There are no upper bounds for $L^{(r)}(1,\chi)$ uniformly for $r\geq 0$ and for all non-principal characters $\chi(\mathrm{mod}\,q)$ with $q\geq 1$ . Recently , Ishikawa [2] has proved that for sufficiently large $r$ , $\left|L^{(r)}(1,\chi)\right|\leq q^{\frac{r}{\log r}-\frac{1}{2}}$ exp $(r\log\ \log r-\frac{r\log\log r}{\log r})$ for a fixed $q\geq 1$ and for a fixed primitive character $\chi(\mathrm{mod}\,q)$ . For $L^{(r)}(1,\chi)$ as a function of $r$ , this is the best result known to date . However , for $L^{(r)}(1,\chi)$ as a function of $q$ , this is a very bad result . In the light of this , our Theorem 3 below shall give a very good upper bound for $\left|L^{(r)}(1,\chi)\right|$ as a function of both $q$ and $r$ simultaneously , though our result is somewhat weaker as a function of $r$ alone . Our Theorem 2 gives simple explicit expressions and upper bounds for Taylor/Laurent coefficients of $\zeta(s,\alpha)$ , $L(s,\chi)$ and $\phi(\lambda,\alpha,s)$ at $s=1$ and 0 . It must be mentioned here that Berndt[1] has , in the case of $\zeta(s,\alpha)$ at $s=1$ , given the following upper bound ,

namely $\left|\gamma_r(\alpha)-(-1)^r \frac{\log^r \alpha}{r!\alpha}\right| \le \frac{4}{r\pi^r}$ for odd $r \ge 1$

$\le \frac{2}{r\pi^r}$ for even $r \ge 2$ .

While we prove

$$\left|\gamma_r(\alpha)-(-1)^r \frac{\log^r \alpha}{r!\alpha}\right| \le e\frac{2^{-r}}{\sqrt{r}} .$$

Next we state our Theorems .

**<u>Theorem 1</u>** : Let $s=\sigma+it$ with $0<\sigma\le 1$ and $t \ge 0$ and let $\psi(u)=u-[u]-\frac{1}{2}$ .

I) For arbitrary $x>0$ , we have

a) $$\zeta(s,\alpha)-\frac{x^{1-s}}{s-1}=\sum_{0\le n\le x-\alpha}(n+\alpha)^{-s}+\psi(x-\alpha)\cdot x^{-s}-s\int_x^\infty \frac{\psi(u-\alpha)}{u^{s+1}}du$$

b) For any integer $r \ge 0$ ,

$$\zeta^{(r)}(s,\alpha)-\frac{d^r}{ds^r}\left(\frac{x^{1-s}}{s-1}\right) = \sum_{0\le n\le x-\alpha}\frac{d^r}{ds^r}(n+\alpha)^{-s}+\psi(x-\alpha)\frac{d^r}{ds^r}(x^{-s})-\int_x^\infty \frac{\psi(u-\alpha)}{u}\frac{d^r}{ds^r}(su^{-s})du .$$

II) Let $a,\ q\ge 1$ be integers with $1\le a\le q$ and let $Z(s,a,q)=q^{-s}\cdot\zeta(s,\frac{a}{q})$ .

Let $X>0$ arbitrary . Then we have

a) $$Z(s,a,q)-\frac{X^{1-s}}{q(s-1)}=\sum_{\substack{n\le X\\ n\equiv a(\mathrm{mod}\,q)}}n^{-s}+\psi(\tfrac{X-a}{q})X^{-s}-s\ q^{-s}\int_{X/q}^\infty \frac{\psi(u-\frac{a}{q})}{u^{s+1}}du$$

b) For any integer $r \ge 0$ ,

$$Z^{(r)}(s,a,q)-\frac{1}{q}\;\frac{d^r}{ds^r}\left(\frac{X^{1-s}}{s-1}\right)=\sum_{\substack{n\le X\\ n\equiv a(\mathrm{mod}\,q)}}\frac{d^r}{ds^r}n^{-s}+\psi\left(\frac{X-a}{q}\right)\frac{d^r}{ds^r}X^{-s}-\int_{X/q}^{\infty}\frac{\psi(u-\frac{a}{q})}{u}\frac{d^r}{ds^r}(s(qu)^{-s})du$$

III) Let $\chi(\mathrm{mod}\,q)$ be a non-principal Dirichlet character and let $X>0$ arbitrary . Then we have

a) $$L(s,\chi)=\sum_{a=1}^{q}\chi(a)\left\{\sum_{\substack{n\le X\\ n\equiv a(\mathrm{mod}\,q)}}n^{-s}+X^{-s}\cdot\psi(\tfrac{X-a}{q})-\int_{X/q}^{\infty}\frac{\psi(u-\frac{a}{q})}{u}(s(qu)^{-s})du\right\}$$

b) For any integer $r\ge 0$,

$$L^{(r)}(s,\chi)=\sum_{a=1}^{q}\chi(a)\left\{\sum_{\substack{n\le X\\ n\equiv a(\mathrm{mod}\,q)}}n^{-s}(-\log n)^r+\psi\left(\frac{X-a}{q}\right)\frac{d^r}{ds^r}X^{-s}-\int_{X/q}^{\infty}\frac{\psi(u-\frac{a}{q})}{u}\frac{d^r}{ds^r}(s(qu)^{-s})du\right\}$$

IV) We have for arbitrary $x>0$ and for $0<\lambda<1$,

a) $\phi(\lambda,\alpha,s)$

$$=\sum_{0\le n\le x-\alpha}e^{2\pi i\lambda n}(n+\alpha)^{-s}+\int_x^{\infty}\frac{e^{2\pi i\lambda(u-\alpha)}}{u^s}du+2\pi i\lambda\int_x^{\infty}\psi(u-\alpha)\frac{e^{2\pi i\lambda(u-\alpha)}}{u^s}du-s\int_x^{\infty}\psi(u-\alpha)\frac{e^{2\pi i\lambda(u-\alpha)}}{u^{s+1}}du$$
$$+\frac{e^{2\pi i\lambda(x-\alpha)}}{x^s}\psi(x-\alpha)$$

b) For any integer $r\ge 0$ , we have

$\phi^{(r)}(\lambda,\alpha,s)$

$$=\sum_{0\le n\le x-\alpha}e^{2\pi i\lambda n}\;((n+\alpha)^{-s})^{(r)}+\int_x^{\infty}e^{2\pi i\lambda(u-\alpha)}(u^{-s})^{(r)}du+2\pi i\lambda\int_x^{\infty}\psi(u-\alpha)e^{2\pi i\lambda(u-\alpha)}\frac{d^r}{ds^r}u^{-s}du$$

$$-\int_x^{\infty}\frac{\psi(u-\alpha)}{u}e^{2\pi i\lambda(u-\alpha)}\frac{d^r}{ds^r}(su^{-s})du+e^{2\pi i\lambda(x-\alpha)}\psi(x-\alpha)\cdot\frac{d^r}{ds^r}x^{-s}\;.$$

**Corollary** : 1) We have for arbitrary $x>0$ ,

$$\lim_{s\to 1}\left(\zeta(s,\alpha)-\tfrac{1}{s-1}\right) = \sum_{0\le n\le x-\alpha}\tfrac{1}{n+\alpha}+\tfrac{\psi(x-\alpha)}{x}-\int_x^\infty \tfrac{\psi(u-\alpha)}{u^2}du \cdot$$

2) For a non-principal character $\chi(\bmod q)$ and for $X>0$ arbitrary ,

$$L(1,\chi)=\sum_{n\le X}\tfrac{\chi(n)}{n}+\sum_{a=1}^{q}\ \chi(a)\left\{\tfrac{\psi(\frac{X-a}{q})}{X}-\tfrac{1}{q}\int_{X/q}^{\infty}\tfrac{\psi(u-\frac{a}{q})}{u^2}du\right\}\ .$$

3) For $x>0$ arbitrary and for $0<\lambda<1$ , we have

$$\phi\ (\lambda,\alpha,1)=\sum_{0\le n\le x-\alpha}\tfrac{e^{2\pi i\lambda n}}{n+\alpha}+\int_x^\infty \tfrac{e^{2\pi i\lambda(u-\alpha)}}{u}du+2\pi i\lambda\int_x^\infty \tfrac{\psi(u-\alpha)}{u}e^{2\pi i\lambda(u-\alpha)}du-\int_x^\infty \tfrac{\psi(u-\alpha)}{u^2}e^{2\pi i\lambda(u-\alpha)}du\ .$$

**<u>Theorem 2</u>** : For any integer $r\ge 1$,

I) we have

a) $$\lim_{s\to 1}\left(\zeta^{(r)}(s,\alpha)-\tfrac{(-1)^r}{(s-1)^{r+1}}\right)=(-1)^r\left\{\tfrac{\log^r\alpha}{\alpha}+r\int_1^\infty \tfrac{\psi(u-\alpha)}{u^2}\log^{r-1}u\,du-\int_1^\infty \tfrac{\psi(u-\alpha)}{u^2}\log^r u\,du\right\}$$

b) If $\zeta(s,\alpha)-\tfrac{1}{s-1}=\sum_{r\ge 0}\gamma_r(\alpha)(s-1)^r$, then

$$\left|\gamma_r(\alpha)-\tfrac{(-1)^r\log^r\alpha}{r!\alpha}\right|\le e\cdot\tfrac{1}{r!}\left(\tfrac{r}{2e}\right)^r\ .$$

II) We have

a) $$\zeta^{(r)}(0,\alpha)+1=(-1)^r\left\{\log^r\alpha+r\int_1^\infty \tfrac{\psi(u-\alpha)}{u}\log^r u\,du\right\}$$

b) If $\zeta(s,\alpha)=\sum_{r\ge 0}\beta_r(\alpha)s^r$ , then

$$\left|\beta_r(\alpha) - \frac{(-1)^r}{r!}\log^r\alpha\right| \le \frac{e}{3}\frac{1}{r!}\left(\frac{r}{e}\right)^r + \frac{1}{r!} \;.$$

III)

a) We have $\phi^{(r)}(\lambda,\alpha,1)$

$$= (-1)^r\left\{\frac{\log^r\alpha}{\alpha} + \int_1^\infty \frac{e^{2\pi i\lambda(u-\alpha)}}{u}\cdot\log^r u du + 2\pi i\lambda\int_1^\infty \frac{\psi(u-\alpha)}{u}e^{2\pi i\lambda(u-\alpha)}\log^r u du - \int_1^\infty \frac{\psi(u-\alpha)}{u^2}e^{2\pi i\lambda(u-\alpha)}\cdot\log^{r-1}u\cdot(r-\log u)du\right\}$$

b) $\left|\frac{\phi^{(r)}(\lambda,\alpha,1)}{r!} - \frac{(-1)^r}{r!}\frac{\log^r\alpha}{\alpha}\right| < \frac{r^r e^{-r}}{r!}\left(\frac{1}{\lambda}+\frac{1}{1-\lambda}\right)$ for $0<\lambda<1$ .

**<u>Theorem</u> 3** : If $\chi(\bmod q)$ is a primitive character modulo $q \ge 3$, then for any integer $r \ge 1$ and for arbitrary $X > 0$ , we have

I)

a) $(-1)^r L^{(r)}(1,\chi)$

$$= \sum_{n\le X}\frac{\chi(n)\log^r n}{n} + \sum_{a=1}^{q}\chi(a)\left\{\psi\left(\frac{X-a}{q}\right)\frac{\log^r X}{X} + \frac{1}{q}\int_{X/q}^\infty \frac{\psi(u-\frac{a}{q})}{u^2}\log^{r-1}qu\cdot(r-\log qu)du\right\}$$

b) $L^{(r)}(1,\chi)$

$$= (-1)^r\sum_{n\le qe^{r/2}}\frac{\chi(n)\log^r n}{n} + 0\left(q^{-1/2}e^{-r/2}\log q\cdot(\log q+\tfrac{r}{2})^r\right) << (\log q+\tfrac{r}{2})^{r+1} \;.$$

II)

a) $(-1)^r L^{(r)}(0,\chi) = \sum_{n\le X}\chi(n)\log^r n + \log^r X\cdot\sum_{a=1}^{q}\chi(a)\left(\frac{X-a}{q} - [\frac{X-a}{q}] - \frac{1}{2}\right)$

$$+ r\sum_{a=1}^{q}\chi(a)\int_{X/q}^\infty \frac{\psi(u-\frac{a}{q})}{u}\log^{r-1}qudu$$

b) $L^{(r)}(0,\chi)=(-1)^r\sum_{n\le qe^{r-1}}\chi(n)\log^r n+0(q^{1/2}\log q\cdot(\log q+r)^r)$

Next we state our Theorem 4 without proof . This Theorem will give complete forms of approximate functional equations of $\zeta(s,\alpha)$ and $L(s,\chi)$ and consequently those of $\zeta^{(r)}(s,\alpha)$ and $L^{(r)}(s,\chi)$ , in the light of our observation below that the expression for $\zeta^{(r)}(s,\alpha)$ or $L^{(r)}(s,\chi)$ can be obtained from that of $\zeta(s,\alpha)$ or $L(s,\chi)$ respectively , by differentiating term-by-term with respect to s , r times or by differentiating under integral sign , if need arises . Note that the complete forms of approximate functional equations for $\zeta$ $(s,\alpha)$ and $L(s,\chi)$ have been stated indirectly in author [4] and author [5] . Our Theorem 4 is useful in the study of $\zeta^{(r)}(s,\alpha)$ as a function of $\alpha$ only , where $0<\alpha\le 1$ ( and consequently , as a function of q in the study of $L^{(r)}(s,\chi)$ ), though not as a function of t unlike the Riemann-Siegel type formula for $\zeta(s,\alpha)$ of author [3] . Before stating our Theorem 4 , it may be noted that the functional equation for $\zeta(s,\alpha)$ can also be written as

$$\zeta(s,\alpha)=2(2\pi)^{s-1}\Gamma(1-s)\sum_{n\ge1}\sin(\tfrac{\pi s}{2}+2\pi n\alpha)n^{s-1}=\Gamma(1-s)\sum_{|n|\ge1}e^{2\pi in\alpha}(2\pi in)^{s-1}\text{ for }\sigma<0.$$

**<u>Theorem</u> 4** : Let $s=\sigma+it$ with $0\le\sigma\le 1$ and $t\ge 0$ .

1) Let $x>0$ arbitrary and let $y=\frac{t}{2\pi x}$ . Then we have

$$\zeta^{(r)}(s,\alpha)-(\tfrac{x^{1-s}}{s-1})^{(r)}=\sum_{0\le n\le x-\alpha}((n+\alpha)^{-s})^{(r)}+\sum_{1\le|n|\le y}e^{2\pi i n\alpha}\Big(\Gamma(1-s)(2\pi i n)^{s-1}\Big)^{(r)}-(x^{-s})^{r}\sum_{|n|>y}\tfrac{e^{2\pi i n(x-\alpha)}}{2\pi i n}$$

$$-\sum_{1\le|n|\le y}\int_0^x e^{2\pi i n(u-\alpha)}\cdot(u^{-s})^{(r)}du+\sum_{|n|>y}\tfrac{1}{2\pi i n}\int_x^\infty \tfrac{e^{2\pi i n(u-\alpha)}}{u}(su^{-s})^{(r)}du$$

2) Let $\chi(\bmod q)$ be a non-principal character modulo an integer $q\ge 1$ and let $X>0$ arbitrary and let $y=\frac{qt}{2\pi X}$ . Then , we have

$$L^{(r)}(s,\chi)=\sum_{a=1}^{q}\chi(a)\Big(Z^{(r)}(s,\tfrac{a}{q})-(\tfrac{X^{1-s}}{q(s-1)})^{(r)}\Big)$$ , where

$$Z^{(r)}(s,\tfrac{a}{q})-(\tfrac{X^{1-s}}{q(s-1)})^{(r)}$$

$$=\sum_{0\le n\le\frac{X-a}{q}}(nq+a)^{-s}+\tfrac{1}{q}\sum_{1\le|n|\le y}e^{\frac{2\pi i n a}{q}}\Big(\Gamma(1-s)(\tfrac{2\pi i n}{q})^{s-1}\Big)^{(r)}+(X^{-s})^{(r)}\sum_{|n|>y}\frac{e^{2\pi i n(\frac{X}{q}-\frac{a}{q})}}{2\pi i n}-\sum_{1\le|n|\le y}\int_0^{\frac{X}{q}}e^{2\pi i n u}\cdot\Big((qu)^{-s}\Big)^{(r)}du$$

$$+\sum_{|n|>y}\int_{X/q}^{\infty}\tfrac{e^{2\pi i n u}}{u}\Big(s(qu)^{-s}\Big)^{(r)}du\ .$$

Next , we give the proofs of our Proposition and our Theorems .

**Proof of Proposition** :

**Method 1** : Using Leibnitz's theorem for differentiation ,

we have $(-1)^r\gamma_r(a,q)=\lim_{s\to1}(q^{-s}(\zeta(s,\frac{a}{q})-\frac{1}{s-1}))^{(r)}=\lim_{s\to1}\sum_{\ell=0}^{r}\binom{r}{\ell}(q^{-s})^{(r-\ell)}\Big(\zeta(s,\frac{a}{q})-\frac{1}{s-1}\Big)^{(\ell)}$ .

$$= \lim_{s\to 1} \sum_{\ell=0}^{r} \binom{r}{\ell} q^{-s}(-\log q)^{r-\ell} \cdot \ell! \gamma_\ell(\tfrac{a}{q}) = \tfrac{1}{q} \sum_{\ell=0}^{r} \binom{r}{\ell} (-1)^{r-\ell} \ell! \gamma_\ell(\tfrac{a}{q}) \log^{r-\ell} q$$

Here the superscript $(\ell)$ denotes $\ell - th$ order derivative with respect to s .

Thus $(-1)^r \gamma_r(a,q) = \tfrac{1}{q} \sum_{\ell=0}^{r} \binom{r}{\ell} (-1)^{r-\ell} \ell!\, \gamma_\ell(\tfrac{a}{q}) \cdot \log^{r-\ell} q$

$$= \tfrac{r!}{q} \sum_{\ell=0}^{r} (-1)^{r-\ell} \tfrac{\log^{r-\ell} q}{(r-\ell)!} \gamma_\ell(\tfrac{a}{q}) = \tfrac{r!}{q} \sum_{k=0}^{r} (-1)^k \tfrac{\log^k q}{k!} \cdot \gamma_{r-k}(\tfrac{a}{q}) \ .$$

**<u>Method</u> 2** : We have $\zeta(s,\alpha) = \tfrac{1}{s-1} + \sum_{n\ge 0} \gamma_n(\alpha)(s-1)^n$

and $q^{-s} = \tfrac{1}{q} e^{-(s-1)\log q} = \tfrac{1}{q} \sum_{n\ge 0} \tfrac{(-1)^n}{n!} \log^n q \cdot (s-1)^n$ .

This gives $q^{-s}\zeta(s,\alpha) = \tfrac{1}{q}\left( \tfrac{1}{s-1} + \sum_{n\ge 0} \gamma_n(\alpha)(s-1)^n \right)\left( \sum_{n\ge 0} \tfrac{(-1)^n}{n!} \log^n q \cdot (s-1)^n \right)$

$$= \tfrac{1}{q} \sum_{n\ge 0} \tfrac{(-1)^n}{n!} \log^n q \cdot (s-1)^{n-1} + \tfrac{1}{q}\left( \sum_{i\ge 0} \gamma_i(\alpha)(s-1)^i \right)\left( \sum_{j\ge 0} \tfrac{(-1)^j}{j!} \log^j q \cdot (s-1)^j \right)$$

$$= \tfrac{1}{q} \sum_{n\ge 0} \tfrac{(-1)^n}{n!} \log^n q \cdot (s-1)^{n-1} + \tfrac{1}{q} \sum_{n\ge 0} c_n(q,\alpha)(s-1)^n ,$$

where $c_n(q,\alpha) = \sum_{j=0}^{n} \gamma_{n-j}(\alpha) \tfrac{(-1)^j}{j!} \log^j q$ .

For a non-principal character $\chi(\bmod q)$,

we have $L(s,\chi) = \sum_{a=1}^{q} \chi(a)\left(q^{-s}\zeta(s,\tfrac{a}{q})\right) = \sum_{a=1}^{q} \chi(a)\left( \tfrac{1}{q} \sum_{n\ge 0} \tfrac{(-1)^n}{n!} \log^n q \cdot (s-1)^{n-1} \right)$

$$+ \sum_{a=1}^{q} \chi(a)\left( \tfrac{1}{q} \sum_{n\ge 0} c_n(q,\tfrac{a}{q})(s-1)^n \right) = \tfrac{1}{q} \sum_{a=1}^{q} \chi(a) \cdot \sum_{n\ge 0} c_n(q,\tfrac{a}{q})(s-1)^n = \sum_{n\ge 0} \gamma_n(\chi)(s-1)^n ,$$

where $\gamma_n(\chi) = \frac{1}{q}\sum_{a=1}^{q}\chi(a)c_n(q,\frac{a}{q}) = \frac{1}{q}\sum_{a=1}^{q}\chi(a)\cdot\sum_{j=0}^{n}\frac{(-1)^j}{j!}\log^j q\cdot\gamma_{n-j}(\frac{a}{q})$

Thus $(-1)^r\gamma_r(a,q) = \frac{r!}{q}\sum_{j=0}^{q}\frac{(-1)^j}{j!}\log^j q\cdot\gamma_{n-j}(\frac{a}{q})$ .

**<u>Proof of Theorem</u> 1 :** We have for $\sigma > 1$, integral $r \geq 0$ and for $0 < \alpha \leq 1$ ,

$$(-1)^r\zeta(s,\alpha) = \sum_{n\geq 0}(n+\alpha)^{-s}\log^r(n+\alpha)$$

$$= \sum_{0\leq n\leq x-\alpha}(n+\alpha)^{-s}\log^r(n+\alpha) + \sum_{n>x-\alpha}(n+\alpha)^{-s}\log^r(n+\alpha)$$ , where $x > 0$ is arbitrary .

Here empty sum has been treated as zero . Also note that $x-\alpha > -1$ .

We apply Euler summation formula to $\sum_{n>x-\alpha}(n+\alpha)^{-s}\log^r(n+\alpha)$ for $\sigma > 1$ .

We have for $\sigma > 1$ , $\sum_{n>x-\alpha}(n+\alpha)^{-s}\log^r(n+\alpha) = \int_{x-\alpha}^{\infty}(u+\alpha)\log^r(u+\alpha)du$

$$+\int_{x-\alpha}^{\infty}\psi(u)\tfrac{d}{du}((u+\alpha)^{-s}\log^r(u+\alpha))du + \tfrac{\psi(x-\alpha)}{x^s}\log^r x \; .$$

Thus for $\sigma > 1$, we have

$$(-1)^r\sum_{n>x-\alpha}(n+\alpha)^{-s}\log^r(n+\alpha) = \int_{x-\alpha}^{\infty}\left(\tfrac{d^r}{ds^r}(u+\alpha)^{-s}\right)du +$$

$$+\int_{x-\alpha}^{\infty}\psi(u)(\tfrac{d}{du}\tfrac{d^r}{ds^r}(u+\alpha)^{-s})du + \int_{x-\alpha}^{\infty}\psi(u)\tfrac{d}{du}(\tfrac{d^r}{ds^r}(u+\alpha)^{-s})du + \psi(x-\alpha)\tfrac{d^r}{ds^r}x^{-s}$$

$$= \tfrac{d^r}{ds^r}\int_{x-\alpha}^{\infty}(u+\alpha)^{-s}du + \int_{x-\alpha}^{\infty}\psi(u)\tfrac{d^r}{ds^r}\left(\tfrac{d}{du}(u+\alpha)^{-s}\right)du + \psi(x-\alpha)\tfrac{d^r}{ds^r}x^{-s}$$

$$= \tfrac{d^r}{ds^r}\left[\tfrac{(u+\alpha)^{1-s}}{1-s}\right]_{u=x-\alpha}^{\infty} + \int_{x-\alpha}^{\infty} \tfrac{d^r}{ds^r}\left(\psi(u)\tfrac{d}{du}(u+\alpha)^{-s}\right)du + \psi(x-\alpha)\cdot\tfrac{d^r}{ds^r}x^{-s}$$

$$= \tfrac{d^r}{ds^r}\left(\tfrac{x^{1-s}}{s-1}\right) + \tfrac{d^r}{ds^r}\left(\int_{x-\alpha}^{\infty}\psi(u)\tfrac{d}{du}(u+\alpha)^{-s}\,du\right) + \tfrac{d^r}{ds^r}\left(\psi(x-\alpha)x^{-s}\right)$$

$$= \tfrac{d^r}{ds^r}\left\{\tfrac{x^{1-s}}{s-1} + \int_{x-\alpha}^{\infty}\psi(u)\tfrac{d}{du}(u+\alpha)^{-s}\,du + \tfrac{\psi(x-\alpha)}{x^s}\right\}$$

$$= \tfrac{d^r}{ds^r}\left(\zeta(s,\alpha) - \sum_{0\le n\le x-\alpha}(n+\alpha)^{-s}\right)$$

However , this is valid for $\sigma > 0$ .

Thus we have for $\sigma > 0$ ,

$$\zeta^{(r)}(s,\alpha) = \sum_{0\le n\le x-\alpha}\tfrac{d^r}{ds^r}(n+\alpha)^{-s} + \tfrac{d^r}{ds^r}\left(\tfrac{x^{1-s}}{s-1}\right) + \tfrac{d^r}{ds^r}\int_{x-\alpha}^{\infty}\psi(u)(-s(u+\alpha)^{-s-1})du + \psi(x-\alpha)\tfrac{d^r}{ds^r}x^{-s}$$

That is, for $\sigma > 0$ we have

$$\zeta^{(r)}(s,\alpha) - \tfrac{d^r}{ds^r}\left(\tfrac{x^{1-s}}{s-1}\right) = \sum_{0\le n\le x-\alpha}\tfrac{d^r}{ds^r}(n+\alpha)^{-s} - \int_{x}^{\infty}\psi(u-\alpha)\tfrac{d^r}{ds^r}(su^{-s-1})du + \psi(x-\alpha)\tfrac{d^r}{ds^r}x^{-s}$$

This proves a) and b) of I) .

II) For given integers a,q with $1 \le a \le q$ and for arbitrary $X > 0$ , this follows from I) on putting $x = \frac{X}{q}$ and $\alpha = \frac{a}{q}$ and by multiplying the expression for $\zeta(s,\frac{a}{q})$ by $q^{-s}$ throughout .

III) This follows from II) on noting $L^{(r)}(s,\chi) = \sum_{a=1}^{q}\chi(a)\cdot Z^{(r)}(s,a,q)$ and on noting

$\sum_{a=1}^{q} \chi(a)\left(\frac{1}{q}\frac{d^r}{ds^r}(\frac{X^{1-s}}{s-1})\right) = 0$ , for any non-principal character $\chi(\bmod q)$ .

IV) Exactly on the same lines as in the case of $\zeta^{(r)}(s,\alpha)$ , we proceed **.**

For $\sigma > 1$, we have $(-1)^r \phi^{(r)}(\lambda,\alpha,s) = \sum_{n\geq 0} e^{2\pi i\lambda n} \cdot (n+\alpha)^{-s} \log^r(n+\alpha)$

$= \left( \sum_{0\leq n\leq x-\alpha} + \sum_{n>x-\alpha} \right) e^{2\pi i\lambda n}(n+\alpha)^{-s} \log^r(n+\alpha)$ , where $x > 0$ is arbitrary .

We apply Euler's summation formula to $\sum_{n>x-\alpha} e^{2\pi i\lambda n}(n+\alpha)^{-s} \log^r(n+\alpha)$ for $\sigma > 1$ and then proceed .

**<u>Proof of Theorem</u> 2** : I) We have for $\sigma > 0$ and for $r \geq 0$ and for $x > 0$ ,

$$\zeta^{(r)}(s,\alpha) - \frac{d^r}{ds^r}\left(\frac{x^{1-s}}{s-1}\right) = \sum_{0\leq n\leq x-\alpha} \frac{d^r}{ds^r}(n+\alpha)^{-s} + \psi(x-\alpha)\frac{d^r}{ds^r}x^{-s} - \int_x^\infty \frac{\psi(u-\alpha)}{u}\frac{d^r}{ds^r}(su^{-s})du \quad .$$

In particular taking $x = 1$ , we have

$$\zeta^{(r)}(s,\alpha) - \frac{d^r}{ds^r}\left(\frac{1}{s-1}\right) = \sum_{0\leq n\leq 1-\alpha} \frac{d^r}{ds^r}(n+\alpha)^{-s} + \psi(1-\alpha)\left(x^{-s}(-\log x)^r\right)_{x=1} - \int_1^\infty \frac{\psi(u-\alpha)}{u}\frac{d^r}{ds^r}(su^{-s})du$$

Thus $\zeta^{(r)}(s,\alpha) - \frac{(-1)^r}{(s-1)^{r+1}} = \frac{d^r}{ds^r}\alpha^{-s} - \int_1^\infty \frac{\psi(u-\alpha)}{u}\frac{d^r}{ds^r}(su^{-s})$ .

Thus $\zeta^{(r)}(s,\alpha) - \frac{(-1)^r}{(s-1)^{r+1}} = \alpha^{-s}(-\log\alpha)^r - \int_1^\infty \frac{\psi(u-\alpha)}{u^{s+1}}(-\log u)^{r-1}(r - s\log u)du$

Letting $s \to 1$ , we get $\lim_{s\to 1}\left(\zeta^{(r)}(s,\alpha) - \frac{(-1)^r}{(s-1)^{r+1}}\right) - \frac{(-\log\alpha)^r}{\alpha} = -\int_1^\infty \frac{\psi(u-\alpha)}{u^2}(-\log u)^{r-1}(r - \log u)du$ .

$$= (-1)^r r \int_1^\infty \frac{\psi(u-\alpha)}{u^2} \log^{r-1} u du - (-1)^r \int_1^\infty \frac{\psi(u-\alpha)}{u^2} \log^r u du \ .$$

Thus $\left\{ \lim_{s\to 1} \left( \zeta^{(r)}(s,\alpha) - \frac{(-1)^r}{(s-1)^{r+1}} \right) - (-1)^r \frac{\log^r \alpha}{\alpha} \right\}$

$$= (-1)^r \left( \int_1^\infty \frac{\psi(u-\alpha)}{u^2} \log^{r-1} u du - \int_1^\infty \frac{\psi(u-\alpha)}{u^2} \log^r u du \right)$$

Next consider $\int_1^\infty \frac{\psi(u-\alpha)}{u^2} \log^{r-1} u du = \left( \int_1^{e^{\frac{r-1}{2}}} + \int_{e^{\frac{r-1}{2}}}^\infty \right) \frac{\psi(u-\alpha)}{u^2} \log^{r-1} u du$ .

Note that $\frac{(\log u)^{r-1}}{u^2}$ is monotonically decreasing for $u > e^{\frac{r-1}{2}}$ and monotonically increasing for $1 < u < e^{\frac{r-1}{2}}$ . This gives , on using second mean value theorem for integrals ,

$$\int_1^{e^{\frac{r-1}{2}}} \frac{\log^{r-1} u}{u^2} \psi(u-\alpha) du = \int_1^{e^{\frac{r-1}{2}}} \frac{\log^{r-1} u}{u^2} d\left( -\sum_{|n|\ge 1} \frac{e^{2\pi i n(u-\alpha)}}{(2\pi i n)^2} \right) = \frac{\log^{r-1}(e^{\frac{r-1}{2}})}{e^{(r-1)}} \int_\xi^{e^{\frac{r-1}{2}}} d\left( \sum_{|n|\ge 1} \frac{e^{2\pi i n(u-\alpha)}}{(2\pi i n)^2} \right)$$

for some $\xi$ with $1 \le \xi \le e^{\frac{r-1}{2}}$ .

Thus $\int_1^{e^{\frac{r-1}{2}}} \frac{\log^{r-1} u}{u^2} \psi(u-\alpha) du \le 4\left(\frac{r-1}{2}\right)^{r-1} \cdot e^{-(r-1)} \sum_{n\ge 1} \frac{1}{4\pi^2 n^2} \le \left(\frac{r-1}{2}\right)^{r-1} \frac{e^{1-r}}{\pi^2} \cdot \frac{\pi^2}{6} \le \frac{e}{6}\left(\frac{r-1}{2}\right)^{r-1} \cdot e^{-r}$ .

Next $\int_{e^{\frac{r-1}{2}}}^\infty \frac{\psi(u-\alpha)}{u^2} \log^{r-1} u du = \lim_{N\to\infty} \int_{e^{\frac{r-1}{2}}}^N \frac{\psi(u-\alpha)}{u^2} \log^{r-1} u du$

$= \lim_{N\to\infty} \int_{e^{\frac{r-1}{2}}}^N \frac{\log^{r-1} u}{u^2} d\psi_2(u-\alpha) du$ , where $\psi_2(u) = -\sum_{|n|\ge 1} \frac{e^{2\pi i n u}}{(2\pi i n)^2}$ .

This gives as before , $\int_{e^{\frac{r-1}{2}}}^\infty \frac{\psi(u-\alpha)}{u^2} \log^{r-1} u du \le \frac{e}{6}\left(\frac{r-1}{2}\right)^{r-1} \cdot e^{-r}$ .

Thus $r\int_1^\infty \frac{\psi(u-\alpha)}{u^2}\log^{r-1} u du \le \frac{er}{3}\left(\frac{r-1}{2}\right)^{r-1} e^{-r} \le \frac{2e}{3}\left(\frac{r}{2}\right)^r e^{-r}$ .

Next , consider $\int_1^\infty \frac{\psi(u-\alpha)}{u^2}\log^r u du = \left(\int_1^{e^{r/2}} + \int_{e^{r/2}}^\infty\right)\frac{\psi(u-\alpha)}{u^2}\cdot\log^r u du$ .

Note that $\frac{\log^r u}{u^2}$ is monotonically decreasing for $u > e^{r/2}$ .

On similar lines as $\int_1^\infty \frac{\psi(u-\alpha)}{u^2}\log^{r-1} u du$ , we find $\int_1^\infty \frac{\psi(u-\alpha)}{u^2}\cdot\log^r u du \le \frac{e}{3}\left(\frac{r}{2}\right)^r e^{-r}$ .

However $\lim_{s\to 1}\left(\zeta^{(r)}(s,\alpha) - \frac{(-1)^r}{(s-1)^{r+1}}\right) = \gamma_r(\alpha)\cdot r!$ .

Thus $\left|r!\gamma_r(\alpha) - (-1)^r\frac{\log^r\alpha}{\alpha}\right| \le e\left(\frac{r}{2}\right)^r e^{-r}$ .

That is , $\left|\gamma_r(\alpha) - (-1)^r\frac{\log^r\alpha}{r!\alpha}\right| \le e\frac{1}{r!}\left(\frac{r}{2e}\right)^r$ .

II) Next we deal with the case $s = 0$ .

From the expression

$$\zeta^{(r)}(s,\alpha) - \frac{(-1)^r}{(s-1)^{r+1}} = \alpha^{-s}(-\log\alpha)^r - \int_1^\infty \frac{\psi(u-\alpha)}{u^{s+1}}(-\log u)^{r-1}(r - s\log u)du ,$$

on putting $s = 0$ , we get $\zeta^{(r)}(0,\alpha) - (-1) = (-\log\alpha)^r - r\int_1^\infty \frac{\psi(u-\alpha)}{u}(-\log u)^{r-1}du$ .

That is , $\zeta^{(r)}(0,\alpha) = -1 + (-1)^r\log^r\alpha + (-1)^r r\int_1^\infty \frac{\psi(u-\alpha)}{u}\log^{r-1} u du$ .

Thus $\zeta^{(r)}(0,\alpha) - (-1)^r\log^r\alpha = -1 + (-1)^r r\left(\int_1^{e^{r-1}} + \int_{e^{r-1}}^\infty\right)\frac{kog^{r-1}u}{u}\cdot\psi(u-\alpha)du$ .

On similar lines as the case $s = 1$,

we have $\left| r\left( \int_1^{e^{r-1}} + \int_{e^{r-1}}^{\infty} \right) \frac{\log^{r-1} u}{u} \psi(u-\alpha)du \right| \le \frac{r}{3}(r-1)^{r-1}e^{-r+1} \le \frac{e}{3} r^r \cdot e^{-r}$

This gives $\left| r!\beta_r(\alpha) - (-1)^r \log^r \alpha \right| \le \frac{e}{3}\left(\frac{r}{e}\right)^r + 1$ .

That is , $\left| \beta_r(\alpha) - \frac{(-1)^r}{r!} \log^r \alpha \right| \le \frac{e}{3}\frac{1}{r!}\left(\frac{r}{e}\right)^r + \frac{1}{r!}$ .

III) We shall deal with the case $\phi(\lambda,\alpha,s)$ at $s=1$ for $0<\lambda<1$ . The case $\lambda=1$ has already been dealt with . For $s>0$ , $0<\lambda<1$ and for $x=1$, we have $\phi^{(r)}(\lambda,\alpha,s) = (\alpha^{-s})^{(r)} + \int_1^{\infty} \left( \frac{e^{2\pi i\lambda(u-\alpha)}}{u^s} \right)^{(r)} du + \int_1^{\infty} \psi(u-\alpha) \frac{d^r}{ds^r} \left( 2\pi i\lambda \cdot \frac{e^{2\pi i\lambda(u-\alpha)}}{u^s} - \frac{se^{2\pi i\lambda(u-\alpha)}}{u^{s+1}} \right) du$

$$= \alpha^{-s}(-\log\alpha)^r + \int_1^{\infty} \left( \frac{e^{2\pi i\lambda(u-\alpha)}}{u^s} \right)(-\log u)^r du + 2\pi i\lambda \int_1^{\infty} \psi(u-\alpha) \frac{e^{2\pi i\lambda(u-\alpha)}}{u^s}(-\log u)^r du$$

$$- \int_1^{\infty} \psi(u-\alpha) \frac{e^{2\pi i\lambda(u-\alpha)}}{u}(su^{-s})^{(r)} du$$

$$= (-1)^r \alpha^{-s} \log^r \alpha + (-1)^r \int_1^{\infty} \frac{e^{2\pi i\lambda(u-\alpha)}}{u^s} \log^r u du + (2\pi i\lambda)(-1)^r \int_1^{\infty} \psi(u-\alpha) \cdot \frac{e^{2\pi i\lambda(u-\alpha)}}{u^s} \log^r u du$$

$$+ \int_1^{\infty} \left( \sum_{|n|\ge 1} \frac{e^{\pi i(n+\lambda)(u-\alpha)}}{2\pi i n} \right) \cdot \frac{(-\log u)^{r-1}}{u^{s+1}} (r - s\log u) du \quad .$$

Thus $\phi^{(r)}(\lambda,\alpha,1) - (-1)^r \frac{\log^r \alpha}{\alpha} = (-1)^r \int_1^{\infty} \frac{e^{2\pi i\lambda(u-\alpha)}}{u} \log^r u du - (-1)^r \int_1^{\infty} \left( \sum_{|n|\ge 1} \frac{\lambda}{n} e^{2\pi i(n+\lambda)(u-\alpha)} \right) \frac{\log^r u}{u} du$

$$- (-1)^r r \int_1^{\infty} \frac{\log^{r-1} u}{u^2} \left( \sum_{|n|\ge 1} \frac{e^{\pi i(n+\lambda)(u-\alpha)}}{2\pi i n} \right) du + (-1)^r \int_1^{\infty} \frac{\log^r u}{u^2} \left( \sum_{|n|\ge 1} \frac{e^{\pi i(n+\lambda)(u-\alpha)}}{2\pi i n} \right) du$$

Thus $(-1)^r \left( \phi^{(r)}(\lambda,\alpha,1) - (-1)^r \frac{\log^r \alpha}{\alpha} \right) = \int_1^{\infty} \frac{\log^r u}{u} d\left( \frac{e^{2\pi i\lambda(u-\alpha)}}{2\pi i\lambda} \right)$

$$-\int_1^\infty \frac{\log^r u}{u}\, d\left(\lambda\cdot\sum_{|n|\geq 1}\frac{e^{2\pi i(n+\lambda)(u-\alpha)}}{2\pi i n(n+\lambda)}\right)-r\int_1^\infty \frac{\log^{r-1} u}{u^2}\cdot d\left(\sum_{|n|\geq 1}\frac{e^{2\pi i(n+\lambda)(u-\alpha)}}{(2\pi i)^2 n(n+\lambda)}\right)+\int_1^\infty \frac{\log^r u}{u^2}\cdot d\left(\sum_{|n|\geq 1}\frac{e^{2\pi i(n+\lambda)(u-\alpha)}}{(2\pi i)^2 n(n+\lambda)}\right)$$

$$=\left(\int_1^{e^r}+\int_{e^r}^\infty\right)\frac{\log^r u}{u}\, d\left(\frac{e^{2\pi i\lambda(u-\alpha)}}{2\pi i\lambda}\right)-\left(\int_1^{e^r}+\int_{e^r}^\infty\right)\frac{\log^r u}{u}\, d\left(\lambda\sum_{|n|\geq 1}\frac{e^{2\pi i(n+\lambda)(u-\alpha)}}{2\pi i n(n+\lambda)}\right)$$

$$-r\left(\int_1^{e^{\frac{r-1}{2}}}+\int_{e^{\frac{r-1}{2}}}^\infty\right)\frac{\log^{r-1} u}{u^2}\cdot d\left(\sum_{|n|\geq 1}\frac{e^{2\pi i(n+\lambda)(u-\alpha)}}{(2\pi i)^2 n(n+\lambda)}\right)+\left(\int_1^{e^{r/2}}+\int_{e^{\frac{r}{2}}}^\infty\right)\frac{\log^r u}{u^2}\, d\left(\sum_{|n|\geq 1}\frac{e^{2\pi i(n+\lambda)(u-\alpha)}}{(2\pi i)^2 n(n+\lambda)}\right)$$

Next , we observe following facts .

a) $\frac{\log^r u}{u}$ is monotonically decreasing for $u > e^r$ and increasing on $[1, e^r]$ .

b) $\frac{\log^{r-1} u}{u^2}$ is monotonically decreasing for $u > e^{\frac{r-1}{2}}$ and increasing on $[1, e^{\frac{r-1}{2}}]$ .

Similarly , $\frac{\log^r u}{u^2}$ is monotonically decreasing for $u > e^{\frac{r}{2}}$ and increasing on $[1, e^{\frac{r}{2}}]$ .

In all the above integrals , we apply second mean value theorem for integrals to the real and imaginary parts of each integral separately . While applying mean value theorem to improper integrals , say to $\int_c^\infty$ for any finite c>0 , we note that $\int_c^\infty = \lim_{N\to\infty}\int_c^N$ .

Consider $\left(\int_1^{e^r}+\int_{e^r}^\infty\right)\frac{\log^r u}{u}\, d\left(\frac{e^{2\pi i\lambda(u-\alpha)}}{2\pi i\lambda}\right)$ .

On using mean value theorem , $\begin{smallmatrix}\mathrm{Re}\\ \mathrm{Im}\end{smallmatrix}\int_1^{e^r}\frac{\log^r u}{u} d\left(\frac{e^{\pi i\lambda(u-\alpha}}{2\pi i\lambda}\right)=\frac{(\log e^r)^r}{e^r}\int_{\xi_1}^{e^r}\begin{smallmatrix}\mathrm{Re}\\ \mathrm{Im}\end{smallmatrix}\; d\left(\frac{e^{2\pi i\lambda(u-\alpha)}}{2\pi i\lambda}\right)$ for some $1\leq \xi_1 \leq e^r$ and this is $=0\left(\frac{r^r e^{-r}}{\lambda}\right)$ .

Next $\lim_{N\to\infty}\int_{e^r}^{N}\frac{\log^r u}{u}d\left(\frac{e^{2\pi i\lambda(u-\alpha)}}{2\pi i\lambda}\right)<<\frac{\left(\log e^r\right)^r}{e^r}\int_{e^r}^{\xi_2}{}^{\mathrm{Re}}_{\mathrm{Im}}d\left(\frac{e^{2\pi i\lambda(u-\alpha)}}{2\pi i\lambda}\right)$ for some real $\xi_2$ and thus $<<\frac{r^r e^{-r}}{\lambda}$ .

Thus $\int_1^\infty\frac{\log^r u}{u}d\left(\frac{e^{2\pi i\lambda(u-\alpha)}}{2\pi i\lambda}\right)<<\frac{r^r e^{-r}}{\lambda}$ .

Similarly , $\lambda\int_1^\infty\frac{\log^r u}{u}d\left(\sum_{|n|\geq 1}\frac{e^{2\pi i(n+\lambda)(u-\alpha)}}{2\pi i n(n+\lambda)}\right)<<r^r e^{-r}\sum_{|n|\geq 1}\frac{1}{n(n+\lambda)}<<\frac{r^r e^{-r}}{1-\lambda}$ .

Similarly , $r\int_1^\infty\frac{\log^{r-1} u}{u^2}d\left(\sum_{|n|\geq 1}\frac{e^{2\pi i(n+\lambda)(u-\alpha)}}{(2\pi i)^2 n(n+\lambda)}\right)<<\frac{r(r-1)^{r-1}}{e^{r-1}}\sum_{|n|\geq 1}\frac{1}{n(n+\lambda)}<<\frac{r^r e^{-r}}{1-\lambda}$ and

$\int_1^\infty\frac{\log^r u}{u^2}d\left(\sum_{|n|\geq 1}\frac{e^{2\pi i(n+\lambda)(u-\alpha)}}{(2\pi i)^2 n(n+\lambda)}\right)<<\frac{r^r e^{-r}}{1-\lambda}$ .

Thus $\left|\phi^{(r)}(\lambda,\alpha,1)-(-1)^r\frac{\log^r\alpha}{\alpha}\right|<<r^r e^{-r}\left(\frac{1}{\lambda}+\frac{1}{1-\lambda}\right)$ .

Hence $\left|\gamma_r(\lambda,\alpha)-\frac{(-1)^r}{r!}\frac{\log^r\alpha}{\alpha}\right|<<\frac{r^r e^{-r}}{r!}\left(\frac{1}{\lambda}+\frac{1}{1-\lambda}\right)$ .

**<u>Proof of Theorem</u> 3** : a) We have for $\sigma>0$ and for $X>0$ arbitrary ,

$$Z^{(r)}(s,\tfrac{a}{q})-\frac{d^r}{ds^r}\left(\frac{X^{1-s}}{q(s-1)}\right)=(-1)^r\sum_{0\leq n\leq\frac{X-a}{q}}(nq+a)^{-s}\log^r(nq+a)+\psi(\tfrac{X-a}{q})\frac{d^r}{ds^r}X^{-s}-\int_{X/q}^{\infty}\frac{\psi(u-\frac{a}{q})}{u}\frac{d^r}{ds^r}(s(qu)^{-s})du$$

Thus for a non-principal character $\chi(\bmod q)$, as

$$L^{(r)}(s,\chi)=\sum_{a=1}^{q}\chi(a)Z^{(r)}(s,\tfrac{a}{q})=\sum_{a=1}^{q}\chi(a)\left(Z^{(r)}(s,\tfrac{a}{q})-\frac{d^r}{ds^r}\left(\frac{X^{1-s}}{q(s-1)}\right)\right)$$

$$=(-1)^r\sum_{n\leq X}\frac{\chi(n)\log^r n}{n^s}+(-1)^r X^{-s}\log^r X\cdot\sum_{a=1}^{q}\chi(a)\psi(\tfrac{X-a}{q})-\sum_{a=1}^{q}\chi(a)\int_{X/q}^{\infty}\frac{\psi(u-\frac{a}{q})}{u}\frac{d^r}{ds^r}\left(s(qu)^{-s}\right)du$$

$$= (-1)^r \sum_{n\le X} \frac{\chi(n)\log^r n}{n^s} + (-1)^r X^{-s}\log^r X \cdot \sum_{a=1}^{q}\chi(a)\left(\tfrac{X-a}{q} - [\tfrac{X-a}{q}] - \tfrac12\right) - \sum_{a=1}^{q}\chi(a)\int_{X/q}^{\infty} \frac{\psi(u-\frac{a}{q})}{u}\frac{(-\log qu)^{r-1}}{(qu)^s}(r - s\log qu)\,du$$

Putting $s=1$, we have $(-1)^r L^{(r)}(1,\chi) = \sum_{n\le X}\frac{\chi(n)\log^r n}{n} + \frac{\log^r X}{X}\sum_{a=1}^{q}\chi(a)\left(\tfrac{X-a}{q} - [\tfrac{X-a}{q}] - \tfrac12\right)$

$$+ \tfrac{1}{q}\sum_{a=1}^{q}\chi(a)\int_{X/q}^{\infty}\frac{\psi(u-\frac{a}{q})}{u^2}\left(\log^{r-1} qu\right)(r - \log qu)\,du \ .$$

We choose $X = qe^{r/2}$ and note that $\frac{\log^r qu}{u^2}$ is monotonically decreasing for

$u > e^{\frac{r}{2}}$ . So also is $\log^{r-1} qu/u^2$ for $u > e^{r/2}$ for $r \ge 1$ .

Now

$$\left|\sum_{n\le qe^{r/2}}\frac{\chi(n)\log^r n}{n}\right| \le \sum_{n\le qe^{r/2}}\frac{\log^r n}{n} \le (\log qe^{r/2})^r \cdot \sum_{n\le qe^{r/2}}\frac{1}{n} \le (\log qe^{r/2})^r(\log qe^{r/2}+1) << (\log q + \tfrac{r}{2})^{r+1}$$

Next with $X = qe^{r/2}$ , we have

$$\frac{\log^r qe^{r/2}}{(qe^{r/2})}\cdot\sum_{a=1}^{q}\chi(a)\left(\tfrac{X-a}{q} - [\tfrac{X-a}{q}] - \tfrac12\right) << \frac{(\log q + \frac{r}{2})^r}{qe^{r/2}}\left(q^{\frac12}\log q\right) << q^{-\frac12}e^{-\frac{r}{2}}\log q\cdot(\log q + \tfrac{r}{2})^r \ .$$

Next , consider $\left|\frac{r}{q}\sum_{a=1}^{q}\chi(a)\int_{e^{r/2}}^{\infty}\psi(u - \tfrac{a}{q})\frac{\log^{r-1} qu}{u^2}\,du\right|$

$$= \left|\frac{r}{q}\int_{e^{r/2}}^{\infty}\sum_{a=1}^{q}\chi(a)\sum_{|n|\ge 1}\frac{e^{2\pi i n(u-\frac{a}{q})}}{2\pi i n}\frac{\log^{r-1} qu}{u^2}\,du\right|$$

$$= \left|\frac{r}{q}\int_{e^{r/2}}^{\infty}\sum_{|n|\ge 1}\frac{e^{2\pi i n u}}{2\pi i n}\left(\sum_{a=1}^{q}\chi(a)e^{-\frac{2\pi i n a}{q}}\right)\frac{\log^{r-1} qu}{u^2}\,du\right|$$

Writing $\tau(\chi,n) = \sum_{a=1}^{q}\chi(a)e^{\frac{2\pi i n a}{q}}$ , we note that $\tau(\chi,n) = \tau(\chi)\cdot\overline{\chi}(n)$ , if $\chi(\bmod q)$ is primitive . Here $\tau(\chi)=\tau(\chi,1)$ .Note that $|\tau(\chi)|=q^{1/2}$ .

Thus $\tau(\chi,-n)=\tau(\chi)\cdot\overline{\chi}(-n)=\chi(-1)\tau(\chi)\cdot\overline{\chi}(n)$

Thus $\left|\frac{r}{q}\int_{e^{r/2}}^{\infty}\left(\sum_{|n\geq 1|}\frac{e^{2\pi inu}\cdot\tau(\chi,-n)}{2\pi in}\right)\frac{\log^{r-1}qu}{u^2}du\right|$

$$=\left|\chi(-1)\tau(\chi)\frac{r}{q}\int_{e^{r/2}}^{\infty}\left(\sum_{|n|\geq 1}\frac{\overline{\chi}(n)e^{2\pi inu}}{2\pi in}\right)\frac{\log^{r-1}qu}{u^2}du\right| << \frac{r}{q}\,|\,\tau(\chi)\,|\int_{e^{r/2}}^{\infty}\frac{\log^{r-1}qu}{u^2}d\left(\sum_{|n|\geq 1}\frac{\overline{\chi}(n)e^{2\pi inu}}{(2\pi in)^2}\right)$$

Using second mean value theorem for integrals , we have the above

$$<<\frac{r}{q}\,|\,\tau(\chi)\,|\cdot\frac{(\log qe^{r/2})^{r-1}}{(e^{r/2})^2} << \frac{r}{q}\cdot q^{\frac{1}{2}}(\log q+\tfrac{r}{2})^{r-1}\cdot e^{-r} << q^{-\frac{1}{2}}\ re^{-r}(\log q+\tfrac{r}{2})^{r-1}$$

Similarly , $\frac{1}{q}\sum_{a=1}^{q}\chi(a)\int_{e^{r/2}}^{\infty}\frac{\log^{r}qu}{u^2}\cdot\psi(u-\frac{a}{q})du << q^{-\frac{1}{2}}e^{-r}\cdot(\log q+\frac{r}{2})^{r}$ .

Thus for primitive character $\chi(\operatorname{mod} q)$ ,where $q\geq 3$ and $r\geq 1$, we have

$$L^{(r)}(1,\chi)=(-1)^r\sum_{n\leq qe^{r/2}}\frac{\chi(n)\log^r n}{n}+0\left(q^{-1/2}\ \log q\cdot e^{-r/2}(\log q+\tfrac{r}{2})^r\right) << (\log q+\tfrac{r}{2})^{r+1}$$

b) Putting $s=0$, we find that

$$(-1)^r L^{(r)}(0,\chi)=\sum_{n\leq X}\chi(n)\log^r n+\log^r X\cdot\sum_{a=1}^{q}\chi(a)\left(\tfrac{X-a}{q}-[\tfrac{X-a}{q}]-\tfrac{1}{2}\right)$$

$+r\sum_{a=1}^{q}\chi(a)\int_{X/q}^{\infty}\frac{\psi(u-\frac{a}{q})}{u}\log^{r-1}qudu$, where $\chi(\operatorname{mod} q)$ is a non-principal character .

Consider $\int_{X/q}^{\infty}\frac{\psi(u-\frac{a}{q})}{u}\log^{r-1}qudu$.

As $\frac{\log^{r-1}qu}{u}$ is monotonically decreasing for $u>e^{r-1}$,

–

if we choose $\frac{X}{q}=e^{r-1}$ , we find that

$$r\sum_{a=1}^{q}\chi(a)\int_{e^{r-1}}^{\infty}\frac{\psi(u-\frac{a}{q})}{u}\log^{r-1} qudu << r\,|\,\tau(\chi)\,|\int_{e^{r-1}}^{\infty}\frac{(\log qe^{r-1})^{r-1}}{e^{r-1}} << q^{\frac{1}{2}}r\cdot(\log q+r-1)^{r-1}\cdot e^{-r}$$

Next $\log^{r} X\cdot\sum_{a=1}^{q}\chi(a)\left(\frac{X-a}{q}-[\frac{X-a}{q}]-\frac{1}{2}\right) << (q^{\frac{1}{2}}\log q)(\log qe^{r-1})^{r} << q^{\frac{1}{2}}\log q\cdot(\log q+r-1)^{r}$

Thus $L^{(r)}(0,\chi)=(-1)^{r}\sum_{n\le qe^{r-1}}\chi(n)\log^{r} n+0\left(q^{1/2}\ \log q\cdot(\log q+r)^{r}\right)$ ,

for any primitive character $\chi(\bmod q)$ .

This completes the proof of Theorem 3 .